\documentclass[a4paper,12pt]{article}

\usepackage{amssymb,latexsym}
\usepackage{euscript}
\usepackage{eufrak}
\usepackage{amsmath}

\newcommand{\negr}[1]{\boldsymbol{#1}}

\newtheorem{coro}{\bf Corollary}[section]
\newtheorem{theo}[coro]{\bf Theorem}

\newtheorem{defi}[coro]{\bf Definition}
\newtheorem{lem}[coro]{\bf Lemma}
\newtheorem{rem}[coro]{\bf Remark}
\newtheorem{ex}[coro]{\bf Example}

\newenvironment{dem}{\noindent\bf Proof. \rm}{\hfill$ \Box$}

\date{}

\title{\Large \textbf{\large \bf  Free algebras in varieties of Hilbert algebras with supremum generated by finite chains. }}

\author{Aldo V. Figallo $^\textup{\scriptsize a}$
\and Elda Pick $^\textup{\scriptsize a}$
\and Susana Saad $^\textup{\scriptsize a}$
 \and Mart\'in Figallo $^\textup{\scriptsize b}$ \footnote{ This work is dedicated to the memory of professor Guillermina Ram\'on who was one of the researchers who initiated the study of these algebras. }
}

\date{
	$^\textup{\scriptsize a}$\textit{\small Instituto de Ciencias B\'asicas, Universidad Nacional de San Juan.}
	\\
	$^\textup{\scriptsize b}$\textit{\small Departamento de Matem\'atica. Universidad Nacional del Sur.}
}

\begin{document}

\maketitle

\thispagestyle{empty}

\begin{abstract}
Hilbert algebras with supremum, i.e., Hilbert algebras where the associated order is a join--semilattice
were first considered by A.V. Figallo, G. Ram\'on and S. Saad in \cite{FRS}, and independently by S. Celani and D.
Montangie in \cite{CeMon}.

On the other hand, L. Monteiro introduced the notion of $n-$valued Hilbert algebras (see \cite{Mon}).
In this work, we investigate the class of $n-$valued Hilbert algebras with supremum, denoted $H^{\vee}_{n}$, i.e.,
$n-$valued Hilbert algebras where the associated order is a join--semilattice. The varieties $H^{\vee}_{n}$ are generated 
by finite chains. The free $H^{\vee}_{n}-$algebra \, ${\bf Free}_{n+1}(r)$\, with $r$ generators is studied. In particular, we determine an upper bound to the cardinal of the finitely generated free algebra \, ${\bf Free}_{n+1}(r)$.

\end{abstract}

\section{\large \bf Introduction}

The study of Hilbert algebras was initiated by Diego in his important work \cite{AD}. It is well--Known that Hilbert algebras constitute an algebraic counterpart of the implicative fragment of Intuitionistic Propositional Logic ({\em IPL}). A topological duality for this algebras was developed in \cite{CeCaMo}.

The class of Hilbert algebras where the associated order is a meet--semilattice was consider in \cite{FavRgSs} under the name of Hilbert algebras with infimum. On the other hand, Hilbert algebras with supremum, i.e., Hilbert algebras where the associated order is a join--semilattice
were first considered by Figallo, Ram\'on and Saad in \cite{FRS}, and independently by Celani and Montangie in \cite{CeMon}. The latter denoted the class of these algebras by $H^{\vee}$.

In 1977, L. Monteiro introduced the notion of $n-$valued Hilbert algebras (see \cite{Mon}). These algebras constitute an algebraic counterpart of the {\em $n-$valued Intuitionistic  Implicative Propositional Calculus}. 

In this work, we investigate the class of $n-$valued Hilbert algebras with supremum (denoted $H^{\vee}_{n}$), i.e., $n-$valued Hilbert algebras where the associated order is a join-semilattice.  The objects of $H^{\vee}_{n}$ are algebraic models for the fragment of intuitionistic propositional calculus in the connectives $\rightarrow$ and $\vee$ and which satisfies the well--known axiom of Ivo Thomas 

$$\beta_{n-1} \rightarrow  (\beta_{n-2} \rightarrow( \dots (\beta_{0}\rightarrow x_{0}) \dots))$$

where \, $\beta_i = (x_{i} \rightarrow x_{i+1})\rightarrow x_{0}, \, 0\leq i \leq n-1$.

It was noted in \cite{CeMon} that $H^{\vee}$ constitute a variety. Here we shall exhibit a very simple and natural equational base for  $H^{\vee}$ different from the one presented in \cite{CeMon}. The most important contribution of this paper is the study of free algebras in $H^{\vee}_{n}$ finitely generated. In particular, we shall describe a formula to determine the cardinal of the free $(n+1)-$valued Hilbert algebra with supremum, \, ${\bf Free}_{n+1}(r)$, \, in terms of the finite number of the free generators $r$ and the numbers $\alpha_{k,p+1}$'s of minimal irreducible deductive systems of certain subalgebras of ${\bf Free}_{n+1}(r)$. Finally, we shall exhibit a formula to calculate an upper bound to $|{\bf Free}_{n+1}(r)|$ in terms of $r$ only.

\section{\large \bf Preliminaries}

If $\mathbb{K}$ is a variety we will denote by $Con_{\mathbb{K}}(A)$, $Hom_{\mathbb{K}}(A,B)$ and $Epi_{\mathbb{K}}(A,B)$ the set of $\mathbb{K}-$congruences of $A$, $\mathbb{K}-$homomorphisms from $A$ into $B$ and  $\mathbb{K}-$epimorphisms from $A$ onto $B$, respectively.

Besides, if $S \subseteq A$ is a $\mathbb{K}-$subalgebra of $A$ we write $S\triangleleft_{\mathbb{K}} A$. We note by $[G]_{\mathbb{K}}$ the $\mathbb{K}-$subalgebra of $A$ generated by the set $G$. When there is no doubt about what variety we are referring to, the subindex will be omitted. 

Recall that a Hilbert algebra is a structure $\langle A, \rightarrow, 1 \rangle$ of type (2,0) that satisfies the following:

\begin{itemize}
\item[(H1)] $x\rightarrow (y\rightarrow x) = 1$,
\item[(H2)] $(x\rightarrow (y\rightarrow z))\rightarrow ((x\rightarrow y) \rightarrow (x\rightarrow z))=1$,
\item[(H3)] $x\rightarrow y = 1 = y \rightarrow x$, implies $x=y$.

\end{itemize}

The variety of all Hilbert algebras is denoted by $H$. It is well-known that in every $A\in H$ the following holds.

\begin{itemize}
\item[(H4)] $(x\rightarrow x)\rightarrow x = x$,
\item[(H5)] $x\rightarrow x =  y\rightarrow y$,
\item[(H6)] $x\rightarrow (y \rightarrow z) = (x\rightarrow y )\rightarrow (x\rightarrow z)$,
\item[(H7)] $(x\rightarrow y)\rightarrow ((y\rightarrow x)\rightarrow x) = (y\rightarrow x)\rightarrow ((x\rightarrow y)\rightarrow y)$.
\item[(H8)] $x\rightarrow (y\rightarrow z) = y\rightarrow (x\rightarrow z)$,
\item[(H9)] $x\rightarrow 1 = 1$,
\end{itemize}

The relation $\leq$ defined by $x \leq y$ iff $x\rightarrow y=1$ is a partial order on $A$ and $x\leq 1$ for all $x \in A$. If $a,b \in A$ are  such that there exists the supremum of \{a,b\} in $A$, denoted by $a\vee b$, then for every $c\in A$ there exists the infimum of \, $\{(a \rightarrow c), (a \rightarrow b) \}$, denoted by  $(a \rightarrow c) \wedge (a \rightarrow b)$  and it is verified that

\begin{itemize} 
\item[(H10)] $ (a\vee b)\rightarrow c = (a\rightarrow c)\wedge (b\rightarrow c).$
\end{itemize}  

If $a, b \in A$  are such that there exists the supremum $a \vee b$ of $\{a, b\}$, then for every $c\in A$ it is verified that:

\begin{itemize} 
\item[(H11)] $(a\rightarrow c)\rightarrow ((b\rightarrow c)\rightarrow ((a \vee b)\rightarrow c))=1.$
\end{itemize}  

It is said that $A \in H$ is a {\em $(n+1)-$valued Hilbert algebra} if $n$ is the least natural number, $n\geq 2$, in such a way that 

\begin{itemize} 
\item[(H12)] $T_{n+1}= T(x_0, \dots, x_n)=\beta_{n-1} \rightarrow (\beta_{n-2} \rightarrow( \dots \rightarrow (\beta_{0} \rightarrow x_{0}) \dots ))=1,$
 
holds in $A$, where $\beta_{i}= (x_{i} \rightarrow x_{i+1})\rightarrow x_{0}$ for $0\leq i \leq n-1$. 
\end{itemize}

We will denote by $H_{n+1}$ the variety of $(n+1)-$valued Hilbert algebras (see \cite{Mon}). 

\begin{ex}\label{Ex1} Let $C_{n+1}=\{0, \frac{1}{n}, \dots, \frac{n-1}{n}, 1 \}$ and let $\rightarrow$ defined as 

$$ x\rightarrow y = \left\{ \begin{tabular}{ll}
$1$ & \mbox{if } $x \leq y$,\\[4mm]
$y$ & \mbox{if } $x > y$.
\end{tabular}\right.$$

Then \, $\langle C_{n+1}, \rightarrow, 1 \rangle$ is a $(n+1)-$valued Hilbert algebra.
\end{ex}

\

Recall that if $A \in H$, $D\subseteq A$ is a {\em deductive system} (d.s.) of $A$ iff \, $1\in D$ \, and $D$ is closed by modus ponens, i.e., if \, $x, x \rightarrow y \in D$ then $y\in D$. A d.s. $D$ is said to be {\em irreducible} (i.d.s.) iff $D$ is a proper d.s. and $D=D_1 \cap D_2 \mbox{ implies } D=D_1 \mbox{ or } D = D_2$ for any d.s. $D_1$ and $D_2$. It is said that the proper d.s. $D$ is {\em fully irreducible} (f.i.d.s) iff 

$$D=\bigcap \limits_{i\in I} D_i \mbox{ implies } D=D_{i} \mbox{ for some } i\in I. (see \cite{AD})$$ 

Besides, $D$ is said to be {\em prime} (p.d.s.) iff it is a proper d.s. and for any $a,b \in A$ we have $a\vee b \in D \mbox{ implies } a \in D \mbox{ or } b \in D$.

We denote by ${\cal D}(A)$ and ${\cal E}(A)$ the set of all d.s and f.i.d.s. of a given Hilbert algebra $A$, respectively. On the other hand, $D$ is {\em minimal} iff it is a minimal element of the ordered set $({\cal D}(A), \subseteq)$, i.e., if $D' \in  {\cal D}(A)$ is such that $D' \subseteq D$ then $D'=D$. By ${\cal M}(A)$ we denote the set of all minimal elements of the ordered set $({\cal E}(A), \subseteq)$, i.e., ${\cal M}(A)=\{D \in {\cal E}(A): \mbox{ if } D' \in  {\cal E}(A) \mbox{ is such that } D' \subseteq D \mbox{ then } D'=D \}$.

The following results are well-known and will be used in the next sections.

\begin{theo} \label{SisDedIr} {\rm( \cite{AD})} Let $A\in H$ and $D\in {\cal D}(A)$. The following conditions are equivalents:
\begin{itemize}
\item[\rm(i)] $D \in {\cal E}(A)$, 
\item[\rm(ii)] if $a,b \in A \setminus D$ then there is $c \in A \setminus D$ such that $a \leq c$ and $b \leq c$,
\end{itemize}  
\end{theo}

\begin{theo} \label{SisDedCompIr}{\rm( \cite{AD})} Let $A\in H$ and $D\in {\cal D}(A)$. The following conditions are equivalents:
\begin{itemize}
\item[\rm(i)] $D$ is a  f.i.d.s., 
\item[\rm(ii)] there is $a\in A \setminus  D$ and $D$ is a maximal d.s. among all d.s. of $A$ that do not contain the element $a$,
\item[\rm (iii)] there is $a\in A \setminus  D$ such that $x \rightarrow a \in D$ for all $x\notin D$.
\end{itemize}  
\end{theo}

Besides,

\begin{theo} \label{MaxSysDed} {\rm(\cite{AD})} Let $A\in H$. For every $D\in {\cal D}(A)$ and $a\notin D$, there exists $M \in {\cal D}(A)$ such that $M$ is maximal among all d.s. that contain $D$ but do not contain $a$.
\end{theo}

The family of all f.i.d.s (minimal i.d.s) of a given Hilbert algebra $A$ is a {\em splitting set}, i.e., 

\begin{equation}\label{EqFamSep}
\bigcap \limits_{D \in {\cal E}(A)} D =\{ 1\} \,\, \, \, \big(\bigcap \limits_{D \in {\cal M}(A)} D =\{ 1\} \big). 
\end{equation}
Recall that if $A, B \in H$ and $h\in Hom_{H}(A,B)$ then  $Con_{H}(A)=\{R(D): D \in {\cal D}(A)\}$ where $R(D)=\{(x,y)\in A^{2}: x\rightarrow y , y \rightarrow x \in D\}$. If $R=R(D)$ we shall denote by $A/D$ the quotient algebra determined by $R$. Also, it is well-known that the {\em kernel} of $h$, $Ker(h)$, is a d.s. of $A$ where $Ker(h)=\{x \in A : h(x)=1\}$. Also, if $z\in A$ the segment $[z)=\{x\in A: z \leq x\}$ is a d.s. of $A$. 

\begin{defi} {\rm(\cite{Mon})} Let $A\in H$. $D\in {\cal D}(A)$ is said to be $(p+1)-$valued if $A/D \simeq C_{p+1}$ (see Example \ref{Ex1}). 
\end{defi}

We shall denote by ${\cal E}_{p+1}(A)$ and  ${\cal M}_{p+1}(A)$ the sets of all $(p+1)-$valued  d.s. of   $A$  included in  ${\cal E}(A)$ and  ${\cal M}(A)$, respectively. Then,

\begin{theo}\label{DSIrNvalued} {\rm (\cite{Mon})} The following conditions are equivalents:
\begin{itemize}
\item[\rm (i)] $A \in H_{n+1}$, 
\item[\rm (ii)] ${\cal E}(A)= \bigcup_{p=1}^{n} {\cal E}_{p+1}(A)$ and ${\cal E}_{n+1}(A)\not= \emptyset$.
\end{itemize}  
\end{theo}

It is worth mentioning that in $H_n$ the notions of i.d.s and f.i.d.s coincide (see \cite{Mon}).
Let $A \in H$ and $X \subseteq A$, we denote by $\mu(X)$ the set of all minimal elements of $X$. Then,

\begin{lem} \label{LemMin}{\rm (\cite{AD})} If $[X]_{H}=A$ then $\mu(X)= \mu(A)$.
\end{lem}

As a consequence of this lemma we have the following corollary.

\begin{coro} Let $A\in H$ and $X\subseteq A$ such that $[X]_{H}=A$. Then, $X=\mu(A)$ iff $X$ is an antichain.
\end{coro}

On the other hand, it can be proved that:

\begin{lem} The only $H$-automorphism of $C_{n+1}$ is the identity.
\end{lem}
\begin{dem} Let $h\in Hom_{H}(C_{n+1}, C_{n+1})$ such that $h$ is a bijection. If $h$ is not the identity map, then there is $x_{0} \in C_{n+1}$ such that $h(x_0)\not= x_0$. Then, $x_0 < h(x_0)$ or $h(x_0)< x_0$. In the first case we have that $x_0 \leq h(x_0) \leq h^{2}(x_0) \leq \dots \leq h^{p}(x_0) \leq \dots$. If there is $p$ such that $h^{p}(x_0)=h^{p+1}(x_0)$ then $x_0=(h^{-1})^{p}(h^{p}(x_0))=(h^{-1})^{p}(h^{p+1}(x_0))=h(x_0)$, a contradiction. So, the only possibility is that $x_0 < h(x_0) < h^{2}(x_0) < \dots < h^{p}(x_0) < \dots $. But this contradicts the finiteness of the algebra $C_{n+1}$.
\end{dem}

\

Finally, we have:

\begin{coro}\label{CorKer} Let $A\in H$ and $h_1, h_2 \in Epi_{H}(A,C_{n+1})$. Then, $h_1=h_2$ iff $Ker(h_1)=Ker(h_2)$.
\end{coro}

\section{\large \bf $\negr{(n+1)-}$valued Hilbert algebras with supremum}

A {\em Hilbert algebra with supremum} (or $H^{\vee}$-algebra), as it was defined in \cite{CeMon}, is an algebra $\langle A, \rightarrow, \vee, 1 \rangle$ of type $(2,2,0)$ such that

\begin{itemize}
\item $\langle A, \rightarrow, 1\rangle$ is a Hilbert algebra,
\item $\langle A, \vee, 1 \rangle$ is a join--semilattice, and
\item for all $a,b \in A$, $a\rightarrow b=1$ iff $a\vee b= b$. 
\end{itemize}

In the next theorem, we exhibit a simple and natural equational base for $H^{\vee}$ different from the one showed in \cite{CeMon}. 

\begin{theo} Let $\langle A, \rightarrow, \vee, 1 \rangle$ be an algebra of type $(2,2,0)$. The following conditions are equivalent.
\begin{itemize}
\item[\rm (i)] $\langle A, \rightarrow, \vee, 1\rangle$ is a $H^{\vee}-$algebra,
\item[\rm (ii)] $\langle A, \rightarrow, 1 \rangle$ is a Hilbert algebra and the following equations hold:\\
{\rm (a)} $x\rightarrow (x\vee y) =1$,\\
{\rm (b)} $y \rightarrow (x\vee y) =1$,\\
{\rm (c)} $(x\rightarrow z) \rightarrow ((y \rightarrow z) \rightarrow ((x \vee y) \rightarrow z))=1$.
\end{itemize}\end{theo}
\begin{dem} It is routine.
\end{dem}

\

Next, we introduce the notion of $(n+1)-$valued Hilbert algebra with supremum. 

\begin{defi}  $A=\langle A, \rightarrow, \vee, 1\rangle\in H^{\vee}$ is a $(n+1)-$valued Hilbert algebra with supremum, or $H^{\vee}_{n+1}-$algebra, if $\langle A, \rightarrow, \vee, 1\rangle \in H^{\vee}$ and $\langle A, \rightarrow, 1\rangle \in H_{n+1}$.
\end{defi}

\begin{ex}\label{Ex2} Let $\langle C_{n+1}, \rightarrow, 1 \rangle$ as in the Example \ref{Ex1} and let $\vee$ the operation defined by $x\vee y =Sup\{x,y\}$. Then,
$J_{n+1}=\langle  C_{n+1}, \rightarrow, \vee, 1 \rangle$ is a $H^{\vee}_{n+1}-$algebra. Besides, if $n \leq m$ then $J_{n}$ is isomorphic to some subalgebra of $J_{m}$.
\end{ex}

It is clear that $H^{\vee}_{n}$ are varieties of Hilbert algebras with supremum generated by finite chains. If fact, $H^{\vee}_{n}$ is generated by the algebra $J_{n}$.
 
All d.s. of a Hilbert algebra with supremum are prime.

\begin{lem}\label{SDCI=SDP} Let $A\in H^{\vee}$ and $D\in {\cal E}(A)$. Then $D$ is prime. 
\end{lem}
\begin{dem} Let $D\in {\cal E}(A)$ and $a,b \in A$ such that (1) $a\vee b \in D$. Suppose that $a\notin A$ and $b\notin A$, by Theorem \ref{SisDedIr}, there exists (2) $c \in A\setminus D$ such that $a \leq c$ and $b \leq c$. Then, \, $1=(a \rightarrow c)\rightarrow((b \rightarrow c) \rightarrow((a \vee b)\rightarrow c))=1 \rightarrow (1 \rightarrow((a \vee b)\rightarrow c)) = (a \vee b)\rightarrow c$. That is, $a \vee b \leq c$ and, by (1), $c \in D$ which contradicts (2).
\end{dem}

\begin{lem}\label{LemHom} Let $h\in Hom_{H^{\vee}}(A,B)$. Then, 
\begin{itemize}
\item[\rm (i)] If $h\in Epi_{H^{\vee}}(A,B)$ then $h([z))=[h(z))$ for all $z \in A$.
\item[\rm (ii)] $Ker(h)\in {\cal E}(A)$, for all $h \in Epi_{H^{\vee}}(A,J_{r+1})$.
\item[\rm (iii)] If $A\in H_{n+1}^{\vee}$ and $B=J_{r+1}$ then $Ker(h)\in {\cal E}_{p+1}(A)$ for some $p$, $1\leq p \leq n$.
\end{itemize}  
\end{lem}
\begin{dem} 

(i) Since $h$ is isotonic we have $h([z))\subseteq [h(z))$. Besides, for all $y \in [h(z))$ there exists $x \in A$ such that $y=h(x)$. Let $u=z \vee x  \in [z)$ then $h(u)=h(z)\vee h(x)=h(z)\vee y= y$. Therefore, $y\in h([z))$.

(ii) Let $a \in A$ such that $h(a)=\frac{r-1}{r}$. Then, $a \notin Ker(h)$ and for all $x\notin Ker(h)$ we have that $h(x\rightarrow a)=h(x)\rightarrow h(a)=h(x) \rightarrow \frac{r-1}{r}= 1 \in Ker(h)$. Then, $x\rightarrow a \in Ker(h)$ and, by Theorem \ref{SisDedCompIr}, $Ker(h)$ is fully irreducible.

(iii) By (ii) and Theorem \ref{DSIrNvalued}.
\end{dem}

\

Theorems \ref{ExtV1} and \ref{ExtV2} extend to $H_{n+1}^{\vee}$ the corresponding results for $H_{n+1}$ proved in \cite{Mon}.

\begin{theo}\label{ExtV1} Let $A\in H_{n+1}^{\vee}$ and $M \in {\cal D}(A)$. The following conditions are equivalent.
\begin{itemize}
\item[\rm (i)] $M\in {\cal E}(A)$,
\item[\rm (ii)] there exists $h \in Hom_{H_{n+1}^{\vee}}(A,J_{r+1})$ such that $Ker(h)=M$, 
\item[\rm (iii)] $A/M$ is isomorphic to some subalgebra of $J_{n+1}$.
\end{itemize}  
\end{theo}

\

Note that if $A\in H^{\vee}$ and $c\in A$ then $[c) \triangleleft_{H^{\vee}} A$. Then,

\begin{lem}\label{DedSys1} Let $A\in H^{\vee}_{n+1}$ and $c\in A$. For every $D \in {\cal E}([c))$ there is a unique $M\in {\cal E}(A)$ such that $D=M\cap [c)$.
\end{lem}
\begin{dem} Let $D \in {\cal E}([c))$. By Theorem \ref{SisDedCompIr}, there is $a \in [c)$ such that $D$ is a maximal d.s. among all d.s. of $A$ that do not contain the element $a$. Then, by Theorem \ref{MaxSysDed}, there exists $M \in {\cal D}(A)$ such that (1) M is maximal among all d.s. of $A$ that verifies:  $a \notin M$ and $D\subseteq M$. Let $P=M \cap [c)$. Then, $D\subseteq P$. 

On the other hand, $P\in{\cal D}([c))$ and $a\notin P$, so, by hypothesis on $D$, $P \subseteq D$. Therefore, $D=P=M\cap[c)$. Besides, suppose that $M=\bigcap \limits_{i\in I} M_i$ ($M_i \in {\cal D}(A)$), then there is $i\in I$ such that $a\notin M_i$. Then, $D \subseteq M \subset M_i$ and, by (1), $M=M_i$. So, $M\in {\cal E}(A)$.

Now, suppose that there are $M_1, M_2 \in {\cal E}(A)$ such that $M_1\cap[c)=M_2\cap[c)=D$. Since $D \not=[c)$, there is (2) $z \in [c)\setminus D$. Then, $z\notin M_1 \cup M_2$. If $z\in M_1\setminus M_2$, by (2), we have $x\vee z \in  M_1\cap [c)\subseteq M_2$. By Lemma \ref{SDCI=SDP} $x\in M_2$ or $z\in M_2$ and both cases lead to a contradiction. Therefore, $M_1 \subseteq M_2$. Analogously, it can be proved that $M_2 \subseteq M_1$. 
\end{dem}

\

In what follows, we shall denote by $M_D$ the only d.s. of $A$ associated to $D \in {\cal E}([c))$.

\begin{lem} Let $A\in H^{\vee}_{n+1}$, $c\in A$, $D\in {\cal E}_{p+1}([c))$ and $M_D\in {\cal E}_{q+1}(A)$. Then, $p\leq q$. Besides, if $M_D\in {\cal M}(A)$ then $D\in {\cal M}([c))$. 
\end{lem}

\begin{dem} By Theorem \ref{ExtV1}, there exists $h \in Hom_{H_{n+1}^{\vee}}(A, J_{n+1})$ such that $Ker(h)=M$. Let $h'=h|_{[c)}$, then $Ker(h')=Ker(h)\cap [c)=M_D\cap [c)=D$ and $p+1=|[c)/D|=|h'([c))|=|h([c))| \leq |h(A)|=|A/M|=q+1$.

Suppose that $M_D\in {\cal M}(A)$ and that there is $D'\in{\cal E}([c))$ such that $D'\subseteq D$. Then, $M_{D'}\cap [c) = D' \subseteq D = M_D \cap [c)$. If there is $x\in M_{D'}\setminus M_D$ then $x\vee z \in M_{D'} \cap [c)=D'$ para cada $z \in [c)\setminus D$, and $x\vee z \in D \subseteq M_D$. By Lemma \ref{SDCI=SDP}, $x \in M_D$ or $z \in M_{D}$ and both cases lead to a contradiction. Then, $M_{D'} \subseteq M_D$ and since $M$ is minimal $M_{D'}=M_D$ and therefore $D'=D$.  
\end{dem}

\

Finally, 

\begin{theo}\label{ExtV2} Let $A\in H_{n+1}^{\vee}$ be a non-trivial algebra. Then, $A$ is isomorphic to a subalgebra of $P=\prod \limits_{M \in {\cal E}(A)} A/M$.
\end{theo}

\

\section{\large \bf Finitely generated free $\negr{H_{n+1}^{\vee}-}$algebras}

Let $r$ be an arbitrary cardinal number, $r>0$. We say that ${\bf Free}(r)$ is the free $H_{n+1}^{\vee}-$algebra with $r$ free generators if:

\begin{itemize}
\item[(L1)] There exists $G \subset {\bf Free}(r)$ such that $|G|=r$ and $[G]_{H_{n+1}^{\vee}} = {\bf Free}(r)$,

\item[(L2)] Any function $f:G\rightarrow A$, $A\in H_{n+1}^{\vee}$, can be extended to a unique homomorphisms $h:{\bf Free}(r) \rightarrow A$.
\end{itemize}

Since the class of $H_{n+1}^{\vee}-$algebras is equationally definable, we know that ${\bf Free}(r)$ is unique up to isomorphisms. If we want to stress that  ${\bf Free}(r) \in H_{n+1}^{\vee}$, we shall write ${\bf Free}_{n+1}(r)$.

\begin{lem} Let $G$ be a set of free generators of ${\bf Free}(r)$. Then, $G=\mu({\bf Free}(r))$.
\end{lem}
\begin{dem} If $G=1$, then $G=\mu(G)$ and, by Lemma \ref{LemMin}, $G=\mu({\cal L}(1))$. Suppose now that $|G|>1$ and let $g, g' \in G$. If $g<g'$, let $f:G\rightarrow J_{n+1}$ the function defined by $f(g)=1$ and $f(t)=0$ if $t\not=g$. By (L2), there exists $h\in Hom_{H_{n+1}^{\vee}}({\bf Free}(r),J_{n+1})$ that extends $f$. Since $h$ is isotonic, $1=h(g)\leq h(g')=0$. Then $g\not<g'$. Analogously, we have that $g'\not< g$. Then, $g$ and $g'$ are incompareble elements and, therefore, $G=\mu(G)=\mu({\bf Free}(r))$. 
\end{dem}

\ 

As an immediate consequence we have:

\begin{coro}\label{CarLib1}  ${\bf Free}(r) = \bigcup \limits_{g \in G} [g)$.
\end{coro}

Taking all this into account we can prove that the variety $H_{n+1}^{\vee}$ is locally finite.

\begin{lem} ${\bf Free}(r)$ is finite, for any natural number $r$.
\end{lem}
\begin{dem} By Theorems \ref{ExtV1} and \ref{ExtV2},it is enough to prove that ${\cal E}({\bf Free}(r))$ is a finite set.
For every $h \in Hom_{H_{n+1}^{\vee}}({\bf Free}(r), J_{n+1})$ we know that $Ker(h)\in {\cal E}({\bf Free}(r))$ and that the correspondence that maps $h$ into $Ker(h)$ is surjective.
Since 

$$Hom({\bf Free}(r), J_{n+1})= \bigcup \limits_{S \triangleleft J_{n+1}} Epi({\bf Free}(r), S),$$

$$\{f=h|_G: h\in Epi({\bf Free}(r), S) \} \subseteq S^{G}$$ 

and $|S^{G}|<\infty$  we have that $|Epi({\bf Free}(r), S)|< \infty$. Since $J_{n+1}$ has a finite number of subalgebras,  ${\bf Free}(r)$ is finite.
\end{dem}

\

Let $G_{k}$ be a subset of $G$ with $k$ elements. Then,

\begin{lem} \label{Lemgk1} Let $g_{k}^{\star} = \bigvee \limits_{g \in G_k} g $. Then, $[g_{k}^{\star})=\bigcap \limits_{g \in G_k} [g)$.
\end{lem}
\begin{dem} It is clear that, if $g\in G_k$ then  $g\leq g_{k}^{\star}$ and $[g_{k}^{\star})\subseteq \bigcap \limits_{g \in G_k} [g)$. 

Let $x \in \bigcap \limits_{g \in G_k} [g)$, then $g\leq x$ for all $g\in G_k$. Then, $g_{k}^{\star} \leq x$ and therefore $x\in [g_{k}^{\star})$. 
\end{dem}

\

\begin{lem} \label{LemImHom} Let  $h\in Epi({\bf Free}(r), J_{q+1})$. Then, $C_{q+1}\setminus \{1\} \subseteq h(G)$.
\end{lem}
\begin{dem} Let $z\in C_{q+1}\setminus \{1\}$ and let $x \in h^{-1}(z) \subseteq {\bf Free}(r)$. Then, we can express $x$ in terms of the generators using the operations $\rightarrow$ and $\vee$. We will call {\em lenght} of $x\in {\bf Free}(r)$ the least natural number $m$ such that there exists an expression for $x$ in terms of the generators which is constructed with $m$ applications of the operations $\rightarrow$ or $\vee$ in it.

If $m=1$, then (1)$x=g_1 \vee g_2$ or (2)$x=g_1 \rightarrow g_2$. In case (1), $h(x)=h(g_1 \vee g_2)= h(g_1) \vee h(g_2)=z$. Since, (3)$h(g_1) \leq h(g_2)$ or $h(g_2) \leq h(g_1)$, we have that $h(g_2)=z$ or $h(g_1)=z$, therefore, $z\in h(G)$. In case (2), if $h(g_1) \leq h(g_2)$ then $h(g_1)\rightarrow h(g_2)=1\not=z$ and so this case is discarded. Then, $h(g_2) \leq h(g_1)$ and if $h(g_1) \rightarrow h(g_2) = h(g_2)=z$ and so $z \in h(G)$.

Suppose that the theorem is true for every formula which its expression in terms of the generators has length $m-1$. Let $x\in {\bf Free}(r)$ be a formula of length $m$. Then,  (4)$x=x_1 \vee x_2$ or (5)$x=x_1 \rightarrow x_2$ with $x_1, x_2\in {\bf Free}(r)$. Clearly, the length of $x_1$and $x_2$ is $m-1$ and, therefore, there exist $g_1, g_2 \in G$ such that $h(g_1)=h(x_1)$ and $h(g_2)=h(x_2)$. In case (4), $h(x)=h(x_1)\vee h(x_2)=h(g_1)\vee h(g_2)=z$ and using the same reasoning above, $h(g_1)=z$ or $h(g_2)=z$. Analogously, in case (5), it must be the case $h(x_2)\leq h(x_1)$ and therefore $h(x_2)=h(x_1)\rightarrow  h(x_2)=z$ and by induction hypothesis there exists $g\in G$ such that $h(g)=h(x_2)=z$.
\end{dem}

\

\begin{lem} Let  $g_{k}^{\star}$ de as in Lemma \ref{Lemgk1}. Then $[g_{k}^{\star})$ is a $H_{m+1}^{\vee}-$subalgebra of ${\bf Free}_{n+1}(r)$ where $m\leq n$.
\end{lem}
\begin{dem} For every $x,y \in [g_{k}^{\star})$, $g_{k}^{\star} \leq y \leq x \rightarrow y$ by (H1), and then $x\rightarrow y \in [g_{k}^{\star})$. Besides, $g_{k}^{\star} \leq y \leq x \vee y$ and so $x\vee y \in [g_{k}^{\star})$. Then, $[g_{k}^{\star}) \triangleleft {\bf Free}_{n+1}(r)$.

Let $m= min\{q: \mbox{ the equation } T_{q+1}\approx 1 \mbox{ holds in } [g_{k}^{\star}) \}$. Then, $m\leq n$ and $[g_{k}^{\star}) \in H_{m+1}^{\vee}$.
\end{dem}

\

Let's denote by $[r]$ the set of the first $r$ natural numbers, i.e., $r=\{1,\dots, r\}$. As an immediate consequence of Corollary \ref{CarLib1} and the {\em inclusion-exclusion principle}  we have:

\begin{lem}\label{CarLibIncExcl} $|{\bf Free}(r)| = \sum \limits_{k=1}^{r} (-1)^{k+1}   \, |\bigcap \limits_{i\in I, \, J\subseteq [r], \, |J|=k} [g_i)|$.
\end{lem}

By Lemma \ref{CarLibIncExcl} and Lemma \ref{Lemgk1} we conclude:

\begin{lem}\label{CarLib2} $|{\bf Free}(r)| = \sum \limits_{k=1}^{r} (-1)^{k+1} \Big(\!\! \begin{array}{l} \scriptstyle{r} \\[-1mm] \scriptstyle{k} \end{array}\!\!\!\Big)   | [g_{k}^{\star})|$.
\end{lem} 

Let ${\cal E}_{k, p+1} = {\cal E}_{p+1}([g_{k}^{\star}))$ and ${\cal M}_{k, p+1}= {\cal M}_{p+1}([g_{k}^{\star}))$.

\begin{rem} \label{Rem1} Let $D\in {\cal E}_{p+1}([g))$ and $h:[g)\longrightarrow J_{p+1}$ be the composition of the canonical map with the isomorphism that there is between $[g)/D$ and $J_{p+1}$. Since $D$ is a $(p+1)-$valued d.s., the family ${\cal S}(D)$ of all d.s of $[g)$ that contain $D$ is a chain whose elements are principal prime d.s. ${\cal S}(D)=\{[a_j): j \in \{0,\dots,p\}\}$. So $D=[a_0) \subseteq [a_1) \subseteq \dots \subseteq [a_j) \subseteq \dots \subseteq [a_p)=[g)$ and $h(a_j)=\frac{p-j}{p}$ for all $j\in\{0, \dots, p\}$ (see \cite{Mon}).
\end{rem}

\

\begin{lem} \label{LemAux1} Let $D\in {\cal M}_{p+1}([g))$, $D'\in {\cal M}_{q+1}([g))$, ${\cal S}(D)=\{[a_j): j \in \{0,\dots,p\}\}$, ${\cal S}(D')=\{[b_i): i \in \{0,\dots,q\}\}$ with $a_0\not= b_0$. Also, let $h_{D'}$ be the canonical epimorphism associated a $D'$ (see Remark \ref{Rem1}). Then,
\begin{itemize}
\item[\rm (i)] $a_j \in [g)\setminus [b_{0})$ for all $j \in \{0,\dots, p\}$,
\item[\rm (ii)] 
$h_{D'}(a_j)= \left\{ \begin{tabular}{ll}
$\frac{p-j}{p} $ & if $D'=D$, \\[4mm]
$0$ & if $D'\not= D$.
\end{tabular}\right.$,

\end{itemize}
\end{lem}
\begin{dem} 

(i) If $a_j \in [b_{0})$, $b_{0} \leq a_j \leq a_0$ and $a_{0} \in [b_{0})$. Then, $[a_0) \subseteq [b_0)$ and since $a_0\not=b_0$, $[b_0)$ is not a minimal f.i.d.s. of $[g)$. So, $a_j \in [g)\setminus [b_{0})$.

(ii) It is  immediate from (i).
\end{dem}

\

\begin{theo}\label{Theogk} $[g_{k}^{\star}) \simeq \prod \limits_{p=1}^{n} J_{p+1}^{\alpha_{k,p+1}}$ \, where \, $\alpha_{k,p+1} = |{\cal M}_{k, p+1}|$.
\end{theo}
\begin{dem} For the sake of simplicity we shall identify $[g_{k}^{\star})/D$ with $J_{p+1}$ for all $D\in {\cal M}_{p+1}([g_{k}^{\star}))= {\cal M}_{p+1}$. Let ${\cal M}= \bigcup \limits_{p=1}^{n} {\cal M}_{p+1}$ and $h_{D} \in Epi_{H_{n+1}^{\vee}}([g_{k}^{\star}), J_{q+1})$ the canonical epimorphism associated to $D$.

Let $h':[g_{k}^{\star})\longrightarrow \prod_{p=1}^{n}  J_{p+1}^{\alpha_{k,p+1}}$ be the function defined as $h'(x)=(h_{D}(x))_{D \in {\cal M}}$. It is clear that $h_{D}$ is a $H_{n+1}^{\vee}-$homomorphism. Besides, $Ker(h')=\{x \in [g_{k}^{\star}):  h'(x)=(h_{D}(x))_{D \in {\cal M}}=1 \}= \bigcap \limits_{D \in {\cal M}} \{x \in [g_{k}^{\star}): h_{D}(x)=1 \} = \bigcap \limits_{D \in {\cal M}} D =\{ 1 \}$, since ${\cal M}$ is a splitting set. Then $h'$ is injective.

Let $y=(y_D)_{D \in {\cal M}} \in \prod \limits_{p=1}^{n} J_{p+1}^{\alpha_{k,p+1}}$. For every $D\in {\cal M}$, $y_D=\frac{p-j}{p}$, for some $j \in \{0, \dots, p\}$. By Remark \ref{Rem1}, there exists $a_D \in [g_{k}^{\star})$ such that $[a_{D})$ is a minimal element in the chain ${\cal S}(D)$ and by Lemma \ref{LemAux1}, $h_{F}(a_{D})=0$ for all $F\in {\cal M}$, $F\not= D$, and $h_D(a_{D})= \frac{p-j}{p}$.

Let $z= \bigvee \limits_{F \in {\cal M}} a_{F}$. Then, $h'(z)=(h_{D}(z))_{D\in {\cal M}} = (h_{D}( \bigvee \limits_{F \in {\cal M}} a_{F}))_{D \in {\cal M}}=
(\bigvee \limits_{F \in {\cal M}} h_{D}( a_{F}))_{D \in {\cal M}}= (h_{D}( a_{D}))_{D \in {\cal M}} = (\frac{p-j}{p})_{D\in {\cal M}}=(y_{D})_{D\in {\cal M}} =y$. So $h'$ is surjective. \end{dem}

\

\begin{coro} Let $G$ a set of free generators. Then,  

$${\bf Free}_{n+1}(r)= \bigcup \limits_{g\in G} [g) \mbox{  where } [g) \simeq \prod \limits_{p=1}^{n} J_{p+1}^{\alpha_{1,p+1}}.$$

\end{coro}

\

\begin{coro} $|{\bf Free}_{n+1}(r)|= \sum \limits_{k=1}^{r} (-1)^{k+1} \Big(\!\! \begin{array}{l} \scriptstyle{r} \\[-1mm] \scriptstyle{k} \end{array}\!\!\!\Big)   \prod \limits_{p=1}^{n} (p+1)^{\alpha_{k,p+1}}$
\end{coro}

\

\section{\large \bf Computing an upper bound to $\negr{\alpha_{k,p+1}}$}

In this section, compute an upper bound to $\alpha_{k,p+1}$. Since ${\cal M}_{k,p+1} \subseteq {\cal E}_{k,p+1}$, $|{\cal M}_{k,p+1}| \leq |{\cal E}_{k,p+1}|$ and so we shall compute the number  $\eta_{k,p+1}=|{\cal E}_{k,p+1}|$. That is to say, we shall determine how many $(p+1)-$valued i.d.s. the subalgebra $[g_{k}^{\star})$ has.

Let $p \leq q \leq n$ and let ${\cal F}_{k,p}(q)$ be the set of all functions $f:G \longrightarrow C_{q+1}$ that satisfies the following conditions:
\begin{itemize}
\item[(F1)]  $C_{q+1} \setminus \{ 1\} \subseteq f(G)$,
\item[(F2)]  $f(g) \leq \frac{q-p}{q}$ \, for all $g\in G_{k}$,
\item[(F3)] there exists $g_p \in G_k$ such that $f(g_p)= \frac{q-p}{q}$.
\end{itemize}

\begin{theo}\label{TeoCardDedSys} $|{\cal E}_{k, p+1}| = \sum \limits_{q=p}^{n} |{\cal F}_{k,p}(q)|$.
\end{theo}
\begin{dem}
Let $f \in \bigcup \limits_{q=p}^{n} {\cal F}_{k,p}(q)$. Then, there exists $q\in \{p, \dots,n \}$ such that $f$ verifies (F1)-(F3). We extend $f$ to a unique $h_f \in Epi_{H_{n+1}^{\vee}}({\bf Free}(r), J_{q+1})$. From (F2) and F(3) we have that $h_f(g_{k}^{\star})= \bigvee \limits_{g\in G_{k}} h_f(g)= \bigvee \limits_{g\in G_{k}} f(g)= \frac{q-p}{q}$. By Lemma \ref{LemHom} (i), $h_f([g_{k}^{\star}))=[h_f(g_{k}^{\star}))=[\frac{q-p}{q})$. Let $D_{f}^{k}=Ker(h_f)\cap [g_{k}^{\star})$. By Theorem \ref{ExtV1}, $Ker(h_f) \in {\cal E}({\bf Free}(r))$. Then, $|[g_{k}^{\star})/D_{f}^{k}|=|h_f([g_{k}^{\star}))|=p+1$ and so  $D_{f}^{k} \in {\cal E}_{k, p+1}$.

Let $\psi: \bigcup \limits_{q=p}^{n} {\cal F}_{k,p}(q) \longrightarrow {\cal E}_{k, p+1}$ be the map defined by $\psi(f)=D_{f}^{k}$. By Theorem \ref{ExtV1}, $\psi$ is well--defined.

For every $D\in {\cal E}_{k, p+1}$, by Lemma \ref{DedSys1}, there is a unique $M\in {\cal E}_{q+1}({\bf Free}(r))$ such that $D=M\cap [g_{k}^{\star})$. Let $h \in Epi_{H_{n+1}^{\vee}}({\bf Free}(r), J_{q+1})$ such that $Ker(h)=M$ (see Theorem \ref{ExtV1}) and let $f$ be the restriction of $h$ to $G$. To see that $f \in {\cal F}_{k,p}(q)$ it is enough to prove that $h(g_{k}^{\star})= \frac{q-p}{q}$.  If $h(g_{k}^{\star})=\frac{i}{q}$ then $p+1=|[g_{k}^{\star})/D|=|[\frac{i}{q})|=q-i+1$. So, $i=q-p$. On the other hand, $h(g_{k}^{\star})=h(\bigvee \limits_{g \in G_k} g)=  \bigvee \limits_{g \in G_k} h(g) = \bigvee \limits_{g \in G_k} f(g)=\frac{q-p}{q}$. Then, it is clear that condition (F2) and (F3) are fulfilled. On the other hand, (F1) holds by Lemma  \ref{LemImHom}. Then, $\psi$ is surjective.

Suppose now that there are $f_1, f_2 \in \bigcup \limits_{q=p}^{n} {\cal F}_{k,p}(q)$ with extensions to ${\bf Free}(r)$ $h_1$ and $h_2$, respectively, such that $D= Ker(h_1)\cap [g_{k}^{\star})=\psi(f_1)=\psi(f_2)= Ker(h_2) \cap [g_{k}^{\star})$. By the uniqueness of the d.s.  $M$ associated to $D$ (Lemma \ref{DedSys1}) we have that $Ker(h_1)=Ker(h_2)$. By  Corollary \ref{CorKer}, $h_1=h_2$ and so $f_1=f_2$. Then, $\psi$ is injective. Taking into account that if $q_1\not= q_2$ then ${\cal F}_{k,p}(q_1) \cap {\cal F}_{k,p}(q_2)=\emptyset$ we conclude that $|{\cal E}_{k, p+1}|=|\sum \limits_{q=p}^{n} {\cal F}_{k,p}(q)| = \sum \limits_{q=p}^{n} |{\cal F}_{k,p}(q)|$.
\end{dem}

\

Next we shall introduce some notation.

\

\begin{equation}
e_{d,a}= \left\{ \begin{tabular}{ll}
$\sum \limits_{j=0}^{a-1} (-1)^{j} \Big(\!\! \begin{array}{c} \scriptstyle{a} \\[-1mm] \scriptstyle{j} \end{array}\!\!\Big) (a-j)^{d}$ & if $1\leq a \leq d$, \\[4mm]
$0$ & if $d<a$ or $a \leq 0$. \end{tabular}\right.
\end{equation}

\

\begin{equation}
u(q,t,b) = e_{r-k,q+1-t+b} + e_{r-k,q-t+b}
\end{equation}

\

\begin{equation}
u(q,t) = \sum \limits_{b=0}^{t} \Big(\!\! \begin{array}{c} \scriptstyle{t} \\[-1mm] \scriptstyle{b} \end{array}\!\!\Big) \, u(q,t,b)
\end{equation}

\

\begin{equation}
\beta(k,p) = \sum \limits_{q=p}^{n} \sum \limits_{t=1}^{q-p+1} \Big(\!\! \begin{array}{c} \scriptstyle{q-p} \\[-1mm] \scriptstyle{t-1} \end{array}\!\!\Big) \, e_{k,t} \cdot u(q,t)
\end{equation}

\

\begin{lem}\label{LemUlt} If $[g_{k}^{\star}) \in H_{m+1}^{\vee}$, $m\leq n$, for all $1\leq k \leq r$, $1\leq p \leq n$, then the number $\eta_{k,p+1}$ of all $(p+1)-$valued i.d.s. of $[g_{k}^{\star})$ verifies:

$$\eta_{k,p+1}= \left\{ \begin{tabular}{ll}
$\beta(k,p)$ & if $p\leq m$,  $k< r$, \\[4mm]
$0$ & if $p>m$,$k< r$,\\[4mm]
$\sum \limits_{q=1}^{r} e_{r,q}$ & if $k=r$,\\[4mm]
 \end{tabular}\right.$$
\end{lem}
\begin{dem} By Theorem \ref{TeoCardDedSys} we know $\eta_{k,p+1}=|{\cal E}_{k, p+1}| = \sum \limits_{q=p}^{n} |{\cal F}_{k,p}(q)|$  for all $1\leq k \leq r$ and $1\leq p \leq n$. 
On the other hand, it is well--known that if $W(A,A')$ is the set of al functions from $A$ into $A'$ is given by

\begin{equation}
|W(A,A')|= \sum \limits_{j=0}^{a-1} (-1)^{j} \Big(\!\! \begin{array}{c} \scriptstyle{a} \\[-1mm] \scriptstyle{j} \end{array}\!\!\Big) (a-j)^{d}
\end{equation}

where $a=|A|$ and $d=|A'|$. 

For $f\in  {\cal F}_{k,p}(q)$ let $f_k = f|_{G_k}$ and $f_{r-k}=f|_{G\setminus G_k}$. By (F2) and (F3), $f_{k}(G_{k})$ verifies that $\frac{q-p}{q} \in f_{k}(G_{k}) \subseteq [0, \frac{q-p}{q}]$. Let ${\cal H}=\{T: T\subseteq [0, \frac{q-p}{q}] \mbox{ such that } \frac{q-p}{q}\in T \}$. Then

\begin{equation}
f_k \in W(G_k, T) \mbox{ for some } T\in {\cal H}
\end{equation}

We want to determine the family of all surjective functions which $f_{r-k}$ belongs to. Let us consider $G_{r-k}=G\setminus G_{k}$, $T=f_{k}(G_{k})$ and $B=f_{r-k}(G_{r-k}) \cap T$. Then, $f_{r-k}(G_{r-k})= \left\{ \begin{tabular}{ll}
$V_{B}=(C_{q+1} \setminus T) \cup B$ & \mbox{ if } $1\in f(G_{r-k})$, \\[4mm]
$V'_{B}=(C_{q+1} \setminus (\{1\}\cup T)) \cup B$ & \mbox{ if } $1\notin f(G_{r-k})$,
\end{tabular}\right.$. So,

\begin{equation}
f_{r-k} \in W(G_{r-k}, V_{B}) \cup W(G_{r-k}, V'_{B})
\end{equation}

for some $B\subseteq T$ and $T\in H$. Let now,

$$ {\cal U}= \bigcup_{T\in {\cal T}} \bigcup_{B \subseteq T} (W(G_{k},T) \times (W(G_{r-k}, V_{B}) \cup W(G_{r-k}, V'_{B}))$$

and

$$\eta :{\cal F}_{k,p}(q) \longrightarrow {\cal U} \mbox{ defined by } \eta(f)=(f_k,f_{r-k}).$$

It is clear that $\eta$ is injective. On the other hand, let $(h_1,h_2) \in {\cal U}$ and consider 
$$h(g)= \left\{ \begin{tabular}{ll}
$h_1(g)$ & \mbox{ if } $g\in G_{k}$, \\[4mm]
$h_2(g)$ & \mbox{ if } $g\in G_{r-k}$.
\end{tabular}\right.$$

If $h_2 \in W(G_{r-k}, V_{B})$ then

\begin{equation}\label{eq1}
h(G)=h_{1}(G_{k}) \cup h_{2}(G_{r-k}) = T \cup V_{B}=C_{q+1}
\end{equation}

If $h_2 \in W(G_{r-k}, V'_{B})$ then

\begin{equation}\label{eq2}
h(G)=h_{1}(G_{k}) \cup h_{2}(G_{r-k}) = T \cup V'_{B}=C_{q+1} \setminus \{1\}
\end{equation}

Besides, if $g\in G_{k}$, $h(g) \in  T \subseteq [0, \frac{q-p}{q}]$ and 

\begin{equation}\label{eq3}
h(g)\leq \frac{q-p}{q} \mbox{ for every } g\in G_{k}.
\end{equation}

Since $h_1(g_p)= \frac{q-p}{q}$ for some $g_{p} \in G_{k}$ it is verified that

\begin{equation}\label{eq4}
h(g_p)=\frac{q-p}{q} \mbox{ for some } g_p\in G_{k}.
\end{equation}

From equations (\ref{eq1}), (\ref{eq2}), (\ref{eq3}) and (\ref{eq4}), we have that $f=h|_{G} \in {\cal F}_{k,p}(q)$ and $\eta$ is surjective.

Then,
\begin{equation}\label{eq5}
|{\cal F}_{k,p}(q)|= |{\cal U}|= |\bigcup_{T\in {\cal T}} \bigcup_{B \subseteq T} (W(G_{k},T) \times (W(G_{r-k}, V_{B}) \cup W(G_{r-k}, V'_{B}))|
\end{equation}

Observe that if $\{A_U \}_{U\in {\cal U}}$  is a family of pairwise disjoint sets such that ${\cal U}=\{U: U\subseteq [0, \frac{m}{q}]\}$, $|A_{U}|=|A_{U'}|$ iff $|U|=|U'|$. If $u=|U|$ then

\begin{equation}
|\bigcup \limits_{U\in {\cal U}} A_{U}|= \sum \limits_{u=0}^{m+1} \Big(\!\! \begin{array}{c} \scriptstyle{m+1} \\[-1mm] \scriptstyle{u} \end{array}\!\!\Big) |A_U|
\end{equation}

Then, if $p\leq m$ and $k<r$

\

\begin{tabular}{lcl}
$\eta_{k,p+1}$ & $=$ & $ \sum \limits_{q=p}^{n} |{\cal F}_{k,p}(q)|$ \\
& $=$ & $ \sum \limits_{q=p}^{n} \sum \limits_{t=1}^{q-p+1} \Big(\!\! \begin{array}{c} \scriptstyle{m+1} \\[-1mm] \scriptstyle{u} \end{array}\!\!\Big) \cdot |W(G_{k},T)| \cdot | (W(G_{r-k}, V_{B}) \cup W(G_{r-k}, V'_{B})|$ \\
& $=$ & $ \sum \limits_{q=p}^{n} \sum \limits_{t=1}^{q-p+1} \Big(\!\! \begin{array}{c} \scriptstyle{m+1} \\[-1mm] \scriptstyle{u} \end{array}\!\!\Big) \cdot  e_{k,t}  \cdot (\sum \limits_{b=0}^{t} \Big(\!\! \begin{array}{c} \scriptstyle{t} \\[-1mm] \scriptstyle{b} \end{array}\!\!\Big) (e_{r-k,q+1-t+b} + e_{r-k, q-t+b}))$\\
& $=$ & $ \sum \limits_{q=p}^{n} \sum \limits_{t=1}^{q-p+1} \Big(\!\! \begin{array}{c} \scriptstyle{m+1} \\[-1mm] \scriptstyle{u} \end{array}\!\!\Big) \cdot e_{k,t}  \cdot (\sum \limits_{b=0}^{t} \Big(\!\! \begin{array}{c} \scriptstyle{t} \\[-1mm] \scriptstyle{b} \end{array}\!\!\Big) u(q,t,b))$\\
& $=$ & $ \sum \limits_{q=p}^{n} \sum \limits_{t=1}^{q-p+1} \Big(\!\! \begin{array}{c} \scriptstyle{m+1} \\[-1mm] \scriptstyle{u} \end{array}\!\!\Big) \cdot e_{k,t}  \cdot u(q,t)$\\
& $=$ & $\beta(k,p)$ \\
\end{tabular}

\

If $k=r$ then $G_{r-k}=\emptyset$ and ${\cal F}_{k,p}(q)=\{ f:G\longrightarrow C_{q+1}: f(G)=[0, \frac{q-p}{q}]=[0, \frac{q-1}{q}]\}$. So, $p=1$ and ${\cal F}_{k,p}(q)=W(G,[\frac{q-1}{q}))$. Then, $\eta_{r,2}=\sum \limits_{q=1}^{r} e_{r,q}$. \end{dem}

\

From Lemma \ref{CarLib2} and Theorem \ref{Theogk} we have

$$|{\bf Free}(r)| \leq \sum \limits_{k=1}^{r} (-1)^{k+1} \Big(\!\! \begin{array}{l} \scriptstyle{r} \\[-1mm] \scriptstyle{k} \end{array}\!\!\!\Big)   \prod \limits_{p=1}^{n} (p+1)^{\eta_{k,p+1}} $$

where $\eta_{k,p+1}$ is as in Lemma \ref{LemUlt}.

\

\begin{ex} 
\begin{itemize}
\item[] 
\item[\rm (i)] For $n=1$ and $r=1$. We have that $\alpha_{1,2}=1$ and
$${\bf Free}_{2}(1) \simeq [g) \simeq J_{2}. $$

\item[\rm (ii)] For $n=1$ and $r=2$, $p=q=1$ and $G=\{g_1,g_2\}$.
$${\bf Free}_{1+1}(2) \simeq [g_1) \cup [g_2) \mbox{ where } [g_1)\simeq [g_2) = J_{2}^{\alpha_{1,2}} = J_{2}^{2}$$

$\alpha_{1,2}=2$ then $[g_i)\simeq J_{2}^{2}$. For $k=2$, $\alpha_{2,2} =1$. On the other hand, $[g_1)\cap [g_2)=[g_1\vee g_2)\simeq J_{2}^{\alpha_{2,2}}=J_2$. Then,

$$|{\bf Free}_{2}(2)| = 6$$
\end{itemize}
\end{ex}

\end{document}